\documentclass[10pt]{article}
\usepackage{amsmath, graphicx}
\usepackage{amssymb}
\usepackage{amsfonts}
\usepackage{amscd}
\usepackage[usenames]{color}
\usepackage{amsthm}
\usepackage{tikz}

\usepackage[top=1.22in, bottom=1.47in, left=1.68in, right=1.68in]{geometry} 
\usepackage{layout}
\usepackage{makecell}
\usepackage{multirow}

\usepackage[hidelinks,draft=false]{hyperref}
\usepackage{bookmark}

\bookmarksetup{
  open,
  numbered,
  addtohook={%
    \ifnum\bookmarkget{level}>2 
      \renewcommand*{\numberline}[1]{}%
    \fi
  },
}

\usepackage{fancybox}
\usepackage{dsfont}
\usepackage{chngcntr}
\counterwithin*{equation}{section}
\usepackage{bm}
\usepackage{mathtools} 
\usepackage{stmaryrd} 
\usepackage{extarrows} 
\usepackage{cancel}

\input{xypic}
\xyoption{all}

\usepackage{threeparttable}
\usepackage{array}
\usepackage{booktabs}
\usepackage{blindtext}
\usepackage{mathrsfs}  

\usepackage{enumitem}
\usepackage{bbm} 
\usepackage{setspace} 
\usepackage{scalerel} 
\usepackage{extsizes} 

\newtheorem{theorem}{Theorem}[section]
\newtheorem{prop}[theorem]{Proposition}

\newtheorem{lemma}[theorem]{Lemma}
\newtheorem{cor}[theorem]{Corollary}

\newtheorem{defn}[theorem]{Definition}

\newtheorem{remark}[theorem]{Remark}
\newenvironment{rem}{\begin{remark}\rm}{\end{remark}}

\newtheorem{facts}[theorem]{Fact}

\newtheorem{example}[theorem]{Example}

\newtheorem{exercise}[theorem]{Exercise}

\newtheorem{terminology}[theorem]{Terminology}

\newtheorem{notation}[theorem]{Notation}

\newtheorem{observation}[theorem]{Observation}

\newtheorem{question}[theorem]{Question}

\makeatletter
\DeclareRobustCommand*\uell{\mathpalette\@uell\relax}
\newcommand*\@uell[2]{
  \setbox0=\hbox{$#1\ell$}
  \setbox1=\hbox{\rotatebox{15}{$#1\ell$}}
  \dimen0=\wd0 \advance\dimen0 by -\wd1 \divide\dimen0 by 2
  \mathord{\lower 0.1ex \hbox{\kern\dimen0\unhbox1\kern\dimen0}}
}

\makeatletter
\renewcommand*\env@matrix[1][*\c@MaxMatrixCols c]{%
  \hskip -\arraycolsep
  \let\@ifnextchar\new@ifnextchar
  \array{#1}}
\makeatother

\def\a{\alpha}
\def\b{\beta}
\def\ga{\gamma}
\def\d{\delta}

\def\s{\sigma}

\def\Ga{\Gamma}
\def\Om{\Omega}

\def\D{\mathscr D}

\def\C{\mathbb C}
\def\P{\mathbb P}

\def\Z{\mathbb Z}
\def\R{\mathbb R}
\def\K{\mathbb K}

\def\lb{\langle}
\def\rb{\rangle}

\def\non{\noindent}
\def\ds{\displaystyle}

\def\pf{\non {\bf Proof. }}

\def\p{\prime}

\def\mqed{\mbox{\qed}}

\def\Hom{{\rm Hom}}

\def\Sym{{\rm Sym}}

\def\pa{\partial}

\newcommand{\sdfrac}[2]{\mbox{\small$\displaystyle\frac{#1}{#2}$}}

\newcommand{\sdbinom}[2]{\mbox{\small$\displaystyle\binom{#1}{#2}$}}

\newcommand{\ip}{\mathbin{\lrcorner}}

\def\S{\mathfrak{S}}

\def\bD{\bar{D}}

\def\bB{\bar{B}}
\def\bK{\bar{K}}
\def\bV{\bar{V}}

\def\H{{\cal H}}
\def\sp{{\rm span}}

\usepackage{scalerel,stackengine}
\newcommand\pig[1]{\scalerel*[5.5pt]{\Big#1}{%
  \ensurestackMath{\addstackgap[1.5pt]{\big#1}}}}

\title{{\bf Witten deformation and divergence-free symmetric Killing 2-tensors }}

\author{Kwangho Choi and   Junho Lee}

\date{\empty}
\addtocounter{section}{0}
\begin{document}

\maketitle

\begin{abstract}

Using a Morse function and a Witten deformation argument, we obtain an upper bound for the dimension of the space of
divergence-free symmetric Killing $p$-tensors on a closed Riemannian manifold, and calculate it explicitly for $p=2$.

\end{abstract}




\section{Introduction}

Let $(M,g)$ be a connected closed Riemannian $n$-manifold and let $\Sym^p(T^*M)$ denote the $p$-th symmetric product of the cotangent bundle $T^*M$. Denote the space of smooth sections of $\Sym^p(T^*M)$ by
$$
S^p(M)=\Ga\big(\Sym^p(T^*M)\big).
$$
Here the symmetric product is given by
$$
\a_1\cdots\a_p=\sum_{\s\in \S_p} \a_{\s(1)}\otimes\cdots\otimes \a_{\s(p)}
\footnote{There is a different convention for the symmetric product $\a\b$ such that a factor $1/p!$ appears on the righthand side (cf. \cite{N}).},
$$
where $\a_1,\cdots,\a_p\in S^1(M)=\Om^1(M)$ and $\S_p$ denotes the symmetric group on $p$ letters.
For simplicity, we will often write $\a\b$ instead of $\a\cdot\b$.

Let $\nabla$ be the Levi-Civita connection on $TM$, and denote its induced connection on $S^p(T^*M)$ by the same $\nabla$. Then the {\em symmetric covariant derivative} is the first order operator
\begin{equation}\label{E:co-v-d}
d^s:S^p(M)\to S^{p+1}(M):\a\mapsto \sum_k e^k \nabla_{e_k}\a,
\end{equation}
where $\{e_k\}$ is an orthonormal frame on $TM$ and its dual coframe $\{e^k\}$.

The Riemannian metric on $M$ induces a pointwise inner product on $S^p(M)$ by the formula
$$
\lb \a_1\cdots\a_p,\b_1\cdots\b_p\rb = \sum_{\s\in \S_p}\lb \a_1,\b_{\s(1)}\rb\cdots\lb \a_p,\b_{\s(p)}\rb,
$$
where $\a_i,\b_j\in S^1(M)$. Denote its induced $L^2$-inner product on $S^p(M)$ by
$$
(\a,\b)=\int_M \lb \a,\b\rb.
$$
Then the formal adjoint of $d^s$ is the {\em divergence operator}
\begin{equation}\label{E:div-op}
d^{s*}:S^{p+1}(M)\to S^p(M):\a\mapsto -\sum_k e_k \ip  \nabla_{e_k}\a
\end{equation}
(cf. Theorem~3.3 of \cite{KP}).

A symmetric tensor $\a\in S^p(M)$ is called {\em Killing tensor} if $d^s\a=0$, and is {\em divergence-free} if $d^{s*}\a=0$. Consider the following self-adjoint second order operator
\begin{equation}\label{E:main-op}
D=d^{s*}d^s + d^sd^{s*}:S^p(M)\to S^p(M).
\end{equation}
It is an elliptic operator on $S^p(M)$ (see Lemma~\ref{L:symbol}) such that
$$
\D^p(M,g):=\ker D = \ker d^s \cap \ker d^{s*},
$$
hence $\D^p(M,g)$ is the space of divergence-free symmetric Killing $p$-tensors.

There is a well-known dimension formula for symmetric Killing $p$-tensors:
\begin{equation}\label{E:DTaTh}
\dim \ker d^s \leq
\frac1n\binom{n+p}{p+1}\sdbinom{n+p-1}{p}
\end{equation}
(see \cite{B,D,MMS,Ta,Th}). The equality holds when $M$ is the $n$-sphere $S^n$ with the standard metric.
When $p=1$, Killing 1-forms are traceless, hence divergence-free (cf. Section~2 of \cite{HMS}).
These are dual to Killing vector fields, and (\ref{E:DTaTh}) specializes to
\begin{equation}\label{E:p=1}
\dim \D^1(M,g)\leq \binom{n+1}{2}.
\end{equation}
This is the well-known bound on the dimension of the Lie algebra of the space of infinitesimal isometries of a Riemannian $n$-manifold.

\bigskip
The aim of this paper is to find a dimension bound, analogous to (\ref{E:DTaTh}), for the space of
{\em divergence-free} symmetric Killing $p$-tensors. Our approach is to use a Morse function $f$ to find
a deformation $D_t$ of the operator $D$
satisfying the following properties:
\begin{itemize}[leftmargin=1.07cm]
\item[(P1)]
$\dim \ker D_t = \dim \ker D$ for all $t$.
\item[(P2)]
$\dim \ker D_t \leq \sum_i \dim\ker \bK(x_i)$ for sufficiently large $t$, where the sum is over the critical
points $\{x_i\}$ of $f$ and $\bK(x_i)$ is a model operator (see Section~\ref{model-operator}).
\item[(P3)]
$\dim\ker \bK(x_i)$ can be calculated.
\end{itemize}
For (P1), we follow Witten's method in \cite{W}; for (P2), we use semiclassical analysis of model operators
in \cite{Sh} (see also \cite{Sh1}); for (P3), we make the technical assumption that {\em the metric $g$ is adapted to $f$}, that is, $g_{ij}=\d_{ij}$ in a Morse chart around each critical point.

Witten deformations and semiclassical analysis are extensively studied and well-understood (see, for example, \cite{CFKS}, \cite{Si}, \cite{Z}, and references therein). However, it is generally difficult to explicitly
calculate the kernel or spectrum of the model operators $\bK(x_i)$. Thus for a deformation $D_t$ satisfying (P1) and (P2), the main task is to calculate $\dim\ker \bK(x_i)$.

\begin{theorem}\label{T:Main}
Let $(M,g)$ be a connected closed $n$-manifold and $f$
be a Morse function such that the metric $g$ is adapted to $f$. Then,
\begin{equation}\label{E:Main-1}
\dim \D^p(M,g)\leq \sum \dim\ker\bK(x_i),
\end{equation}
where the sum is over all critical points $x_i$ of $f$.
Furthermore, for the critical point $x$ of $f$ with index $m$,  the integers
\begin{equation*}\label{E:Main-2}
K^{n,p}_m=\dim\ker\bK(x)
\end{equation*}
depend only on $n$, $p$, and $m$, and  satisfy
$$
\begin{array}{l}
(a)\
K^{n,p}_m=0\ \text{if $m\geq p$ or if $m=n$}.
\\[10pt]
(b)\
K^{n,1}_0=  \sdbinom{n+1}{2}\ \ \
(c)\  K^{n,2}_0 ={\small \left\{
\begin{array}{cl}
1 &\text{if}\ n=1
\\
3 &\text{if}\ n=2
\\
\sdfrac1n\sdbinom{n+2}{3}\sdbinom{n+1}{2}-\sdfrac{n(n+3)}{2} &\text{if}\ n\geq 3
\end{array}\right.}
\end{array}
$$
\end{theorem}

The inequality ({\ref{E:Main-1}) is a consequence of  (P1) and (P2), which are proved in Proposition~\ref{P:WD}\,(a)
and Proposition~\ref{P:Main}, respectively, and the integers $K^{n,p}_m$ depend only on $n$, $p$, and $m$ by
(\ref{E:MO-1}).
When $p=1$, this theorem implies the inequality (\ref{E:p=1}) under the adapted metric assumption because every connected closed manifold admits a Morse function with a unique critical point of index zero.
When $p=2$, Theorem~\ref{T:Main} shows:

\begin{cor}\label{C:Main}
If $M$ admits a Morse function with a unique minimum and no critical point of index $1$ then, for an adapted metric,
$$
\dim\D^2(M,g)\leq K^{n,2}_0,
$$
where $K^{n,2}_0$ is given above.
\end{cor}

A simply connected closed manifold admitting a perfect Morse function satisfies the condition of Corollary~\ref{C:Main}. Examples of such a manifold
include $S^n$ ($n\geq 2$), $\C\P^n$, complex Grassmannian manifolds and their products (cf. \cite{N1}).

\smallskip
At the moment, the authors do not know whether the adapted metric assumption can be eliminated and
whether $K^{n,p}_m=0$ when $m\ne 0$

\medskip
The rest of the paper is organized as follows.
Section~2 reviews basic operators on $\Sym^p(TM)$ and calculates the symbol of $D$. In Section~3, we define a deformation $D_t$ satisfying (P1) and find its model operators $\bK(x)$ satisfying (P2). We then employ Hermite polynomials and operators satisfying the canonical commutation relations to prove Theorem~\ref{T:Main}\,(a) in Section~4, and to compute $K^{n,1}_0$ and $K^{n,2}_0$ in Section~5; this completes the proof of Theorem~\ref{T:Main}, hence shows (P3).

\medskip
\non
{\bf Acknowledgement:} The authors sincerely thank Thomas H. Parker for his insightful discussion and helpful comments on this paper. The second author also thanks Seoul National University and the Korea Institute for Advanced Study for their support and hospitality during his visit when work on this paper was undertaken.

\section{Preliminaries}

This section reviews the canonical commutation relations (\ref{E:CCR}) below, fixes notations, and calculates the symbol of the operator
$D$ defined in (\ref{E:main-op}).

\label{CCR}

\subsection{Canonical commutation relations}

\begin{defn}
For each 1-form $\ga\in S^1(M)=\Om^1(M)$, define operators
\begin{align*}
c_\ga:S^p(M)&\to S^{p+1}(M)\ \ \ \text{by}\ \ \ c_\ga(\a) =\ga\a,
\\
a_\ga:S^p(M)&\to S^{p-1}(M)\ \ \ \text{by}\ \ \ a_\ga(\a) = \ga^*\!\ip\a,
\end{align*}
where $\ga^*$ is the metric dual of $\ga$, i.e., $\ga(v)=g(\ga^*,v)$.
\end{defn}

Let $\{e^k\}$ be an orthonormal frame on $T^*M$.
If we regard $S^p(M)$ as the space of polynomials of degree $p$ in the variables $e^1,\cdots,e^n$, then $c_\ga$ is the multiplication by $\ga$ and $a_\ga$ is the directional derivative in the direction of $\ga$. For instance, if $\ga=\sum_k \ga_ke^k$, $\a=(e^1)^3e^2$, and $\b=(e^1)^2e^2$ then
$$
\lb a_\ga\a,\b\rb =\lb 3\ga_1(e^1)^2e^2+ \ga_2(e^1)^3,\b\rb=6\ga_1=\lb \a,c_\ga\b\rb.
$$

The following well-known fact plays an essential role in our discussions.

\begin{lemma}\label{L:CCR}
For any $\ga,\mu\in S^1(M)$, we have:
\begin{itemize}
\item[(a)]
The adjoint of $a_\ga$ is $a_\ga^*=c_\ga$.
\item[(b)]
The operators $c_\ga$ and $a_\ga$ satisfy the canonical commutation relations
\begin{equation}\label{E:CCR}
c_\ga c_\mu=c_\mu c_\ga,\ \ \ \ a_\ga a_\mu =a_\mu a_\ga,\ \ \ \
c_\ga a_\mu -a_\mu c_\ga =-\lb \ga,\mu \rb Id.
\end{equation}
\end{itemize}
\end{lemma}

Let $\nabla$ denote the connection on $\Sym^p(T^*M)$ induced from the Levi-Civita connection.
Given an orthonormal frame $\{e_k\}$ on $TM$ and its dual coframe $\{e^k\}$, we set
$$
\nabla_k=\nabla_{e_k},\ \ \ \ c_k=c_{e^k},\ \ \ \ a_k=a_{e^k}.
$$
Then, we can write the symmetric covariant derivative $d^s$ in (\ref{E:co-v-d}) and the divergence operator $d^{s*}$
in (\ref{E:div-op}) as
$$
d^s = \sum_k c_k\nabla_k\ \ \ \ \text{and}\ \ \ \
d^{s*}= - \sum_k a_k\nabla_{k}.
$$

\subsection{The symbol map}

\label{SM}

Let  $E$ and $F$ be vector bundles over $M$.
Recall that the (principal) {\em symbol} of a $k^{th}$-order liner differential operator $D:\Ga(E)\to \Ga(F)$
is the section
$$
\s_D\in\Ga\big(\Hom(\pi^*E,\pi^*F)\big)\ \ \ \ \text{(where $\pi:T^*M\to M$)}
$$
such that for $(x,w)\in T^*M$, the symbol map $\s_D(x,w):E_x\to F_x$ is given by
$$
\s_D(x,w)(\xi)=\sdfrac{1}{k!}D(f^k\xi)(x),
\footnote{Some authors use a different convention in which a factor $(-i)^k$ appears on the right-hand side (cf. \cite{BGV})}
$$
where $f$ is a function satisfying $f(x)=0$ and  $df(x)=w$.
A linear differential operator $D$ is called {\em elliptic}
if the symbol map $\s_D(x,w)$ is an isomorphism for every $(x,w)\in T^*M\setminus \{0\}$.

\begin{lemma}\label{L:symbol}
\mbox{}
\begin{itemize}
\item[(a)]
$d^{*s}d^s$ is elliptic, hence $\ker d^s=\ker d^{s*}d^s$ has a finite dimension.
\item[(b)]
$d^{s*}d^s \pm d^{s}d^{s*}$ are elliptic.
\end{itemize}

\end{lemma}

\pf
The symbols of $d^s$, $d^{s*}$, $d^{s*}d^s$, and $d^sd^{s*}$ are given as
\begin{equation}\label{E:symbol}
\begin{array}{ll}
\s_{d^s}(x,w)=c_w,  &\s_{d^{s^*}}(x,w)=-a_w,
\\
\s_{d^{s*}d^s}(x,w)=-a_w c_w,\ \ \  &\s_{d^sd^{s*}}(x,w)=-c_w a_w,
\end{array}
\end{equation}
where $(x,w)\in T^*M\setminus \{0\}$. Since $c_w$ is injective and $a_w^*=c_w$, the symbol $-a_w c_w$ is also injective, hence an isomorphism. This shows (a). For (b), first note that
$$
c_wa_w-a_wc_w=-|w|^2Id
$$
by (\ref{E:CCR}). Also, note that $c_wa_w+a_wc_w$ is an isomorphism because
$$
\big\lb (a_wc_w+c_wa_w)\a,\a\big \rb=\lb c_w\a,c_w\a\rb + \lb a_w\a,a_w\a\rb
$$
and $c_w$ is injective. Therefore, (b) follows.
 \mqed

\begin{rem}
In general, the divergence operator $d^{s*}$ has an infinite-dimensional kernel since its symbol is given by contraction
$a_w$. One can find such an example in Section~3.3 of \cite{CGGR}
in the context of symplectic Dirac operators.
\end{rem}

\begin{rem}
The elliptic operator $d^{s*}d^s - d^{s}d^{s*}$ on symmetric tensors, first studied by Sampson \cite{S}, is a zeroth-order perturbation of the trace Laplacian; hence it is a symmetric analog of the Hodge Laplacian on differential forms.
\end{rem}

\section{Witten deformation}


This section applies Witten's method in \cite{W} to perturb the operator $D=d^{s*}d^s + d^sd^{s*}$ in (\ref{E:main-op}) via a Morse function $f$ and then uses Corollary~1.3 of \cite{Sh} to obtain an upper bound for $\dim \ker D$ given by local data (or model operators) at the critical points of $f$.

\subsection{Perturbation}

We will choose a deformation $D_t$ of $D$ satisfying (P1) and (P2) in the Introduction as follows.
Fix a smooth function $f$ on $M$ and for $t\in\R$, define
$$
d_{t}:S^p(M)\to S^{p+1}(M)\ \ \ \text{by}\ \ \ d_t=e^{-tf}d^se^{tf}.
$$

\begin{defn}
For $t\in\R$, let $D_t= d_{t}^*d_t + d_{-t}d_{-t}^*$.
\end{defn}

Since $d_t= d^s + tc_{df}$ and $d_t^*=e^{tf}d^{s*}e^{-tf}=d^{s*} + ta_{df}$,
we have
\begin{equation}\label{E:PO}
\begin{aligned}
&D_t= D+ tB+t^2V, \ \ \ \text{where}
\\[5pt]
&B=d^{s*} c_{df}-c_{df}d^{s*} +a_{df}d^s- d^s a_{df},\ \ \
V=a_{df}c_{df}+c_{df}a_{df}.
\end{aligned}
\end{equation}
Note that $D_t$, $D$, $B$, and $V$ are all self-adjoint.
The proposition below shows (P1) and is crucial for obtaining (P2).

\begin{prop}\label{P:WD}
\mbox{}
\begin{itemize}
\item[(a)]
$\dim\ker D_t=\dim\ker D$ for all $t\in \R$.
\item[(b)]
$B$ is a bundle map, i.e., its symbol $\s_B=0$.
\item[(c)]
$V=2a_{df}^*a_{df}+|df|^2Id$, hence $\lb V\a,\a\rb \geq |df|^2|\a|^2$, $\forall \a\in S^p(M)$.

\end{itemize}
\end{prop}

\pf
Define an isomorphism $L_{t}:S^p(M)\to S^p(M)$ by
$$
L_t(\a)=e^{-tf}\a.
$$
Then, for $d^s,d_t:S^p(M)\to S^{p+1}(M)$ and $d^{s*},d_{-t}^*:S^p(M)\to S^{p-1}(M)$,
$$
L_{t}(\ker d^s)=\ker d_t\ \ \ \ \text{and}\ \ \ \ L_{t}(\ker d^{s*})=\ker d^*_{-t}.
$$
Now, (a) follows from
$$
\ker D_t = \ker d_t\cap \ker d_{-t}^* = L_t\big(\ker d^s\cap \ker d^{s*}\big)=L_t\big(\ker D\big).
$$
(b) follows from (\ref{E:symbol}) and (\ref{E:CCR}):
$$
\s_B(x,w)=-a_w c_{df}+c_{df}a_w +a_{df}c_w-c_w a_{df}= -\lb df,w\rb Id + \lb w, df\rb Id =0.
$$
Lastly, (c) follows from Lemma~\ref{L:CCR}. \mqed

\medskip

Because $B$ is a bundle map, $D_t$ is a zeroth-order perturbation of $D$ such that,
as $\lb D\a,\a\rb\geq 0$ and $\lb V\a,\a\rb\geq |df|^2|\a|^2$,
elements of $\ker D_t$ concentrate at critical points of $f$ as $t\to\infty$
(cf. Lemma~9.2 of \cite{LP} or Proposition~2.4 of \cite{M}).
For further examples of concentration by perturbing elliptic operators, see Section~7 of \cite{T}, Section~9 of \cite{LP}, \cite{PR}, and \cite{M}.

\subsection{Local data}

In this subsection, we compute model operators under the adapted metric assumption below
and apply Corollary~1.3 of \cite{Sh} to show (P2).

\subsubsection{Adapted metric}
\label{Adapted metric}

Suppose $f$ is a Morse function and $x$ is a critical point of index $m$. Then, by Morse Lemma,
there exists local coordinates $\{y_k\}$ around $x$ with $y_k(x)=0$ such that
$$
f(y)=f(0)-\sdfrac12 y_1^2-\cdots-\sdfrac12 y_m^2 + \sdfrac12 y_{m+1}^2+\cdots +\sdfrac12 y_n^2.
$$
Now assume that the Riemannian metric $g$ on $M$ is {\em adapted to $f$}, that is,
the metric $g$ is given by
$$
g = (dy_1)^2 +\cdots + (dy_n)^2.
$$
(cf. \cite{N1}). Such a metric always exists and the 1-forms $e^k=dy_k$ form an orthonormal frame of $T^*M$.

Let $\a=\sum_J \a_J e^J$, where $J=(j_1,\cdots,j_n)$ is a multi-index, that is,
$e^J=(e^1)^{j_1}\cdots(e^n)^{j_n}.$
Then, in the local coordinates $\{y_k\}$,
\begin{equation}\label{E:simple-cal}
\begin{aligned}
&d^s\a=\sum_J\sum_k \sdfrac{\pa \a_J}{\pa y_k} c_k e^J,\ \ \ \
d^{s*}\a=-\sum_J\sum_k \sdfrac{\pa \a_J}{\pa y_k} a_k e^J,
\\
&df=-y^1e^1\cdots-y^{m}e^m+y^{m+1}e^{m+1}+\cdots+y^ne^n,
\\
&\nabla_k df = \sum_k s_ke^k\ \ \ \text{where}\ \ \
s_k=\left\{
\begin{array}{rl}
-1\ \ &k\leq m,
\\
1\ \  &k>m
\end{array}
\right.
\\
&V(y)=\sum_{k,\ell}(s_ks_\ell)y_ky_\ell(a_kc_\ell+c_ka_\ell).
\end{aligned}
\end{equation}
Observe that by Proposition~\ref{P:WD}\,(c), we have
\begin{equation}\label{E:(C)}
V(0)=0\ \ \ \ \text{and}\ \ \ \ \lb V(y)\a,\a\rb\geq |y|^2|\a|^2.
\end{equation}

\subsubsection{Model operators}

\label{model-operator}

Recall that the rank of $\Sym^p(T^*M)$ is the dimension of the space of homogeneous polynomials of degree
$p$ on $\R^n$, which is
\begin{equation}\label{E:dim}
n_p:=\binom{n+p-1}{p}.
\end{equation}
In the local coordinates $\{y_k\}$ and a trivialization of $\Sym^p(T^*M)$ given by the local frame $\{e^k=dy_k\}$,
the operator $D_t=D+tB+t^2V$ becomes an $n_p\times n_p$ matrix differential operator in $\R^n$ and
$D$ becomes a matrix differential operator of the form
$$
D=\sum_{k,\ell}\frac{\pa^2}{\pa y_k\pa y_\ell} D_{k\ell}(y)  + \text{lower order terms}.
$$
The model operator $\bK(x)$ at the critical point $x$
is  the sum of operators
\begin{equation*}\label{E:model}
\bK(x)=\bD+\bB+\bV,
\end{equation*}
where  the operators $\bD$, $\bB$, and $\bV$ are defined as:
\begin{itemize}
\item
${\ds \bD=\sum_{k,\ell}\sdfrac{\pa^2}{\pa y_k\pa y_\ell} D_{k\ell}(0)}$ is the second order homogeneous matrix differential operator obtained by taking the second order terms of $D$ with the coefficients frozen at the point $x$.

\item
$\bB=B(0)$ is an endomorphism of the fiber $\Sym^p(T^*_xM)$.
\item
${\ds \bV=\sdfrac12\sum_{k,\ell}\sdfrac{\pa^2 V}{\pa y_k\pa y_\ell}(0)y_ky_\ell}$ is the quadratic part of $V$ near $x$.
\end{itemize}
(See \cite{Sh} and also \cite{Sh1}.) In our case, by (\ref{E:CCR}) and
(\ref{E:simple-cal}),
\begin{equation}\label{E:MO-1}
\begin{aligned}
\bar{D} &=-\sum_{k,\ell} \pa_k\pa_\ell(a_kc_\ell+c_ka_\ell),\ \ \ \
\bar{B} 
= -\sum_{k} s_k(a_kc_k+c_ka_k),
\\
\bar{V} 
&=\sum_{k,\ell}(s_ks_\ell)y_ky_\ell(a_kc_\ell+c_ka_\ell).
\end{aligned}
\end{equation}

The model operator $\bK(x)$ is a self-adjoint second order elliptic operator acting in $L^2(\R^n)\otimes \Sym^p(T^*_xM)$, as a closure
from the domain $C_0^\infty(\R^n)\otimes \Sym^p(T^*_xM)$.
Moreover, $\bK(x)$ is positive semi-definite (see Remark~\ref{R:unweighted} below).
Observe that (\ref{E:MO-1}) implies the dimension $K^{n,p}_m=\dim \ker\bK(x)$ depends only on $n$, $p$, and $m$, independent of the critical point $x$.

\subsubsection{Semiclassical analysis}

Let $x_1,\cdots, x_N$ be the critical points of the Morse function $f$. Then the model operator for $D_t$ on $M$ is the elliptic operator
$$
K=\bigoplus_{i=1}^N \bK(x_i)
$$
acting in the orthogonal direct sum $\bigoplus_i \big(L^2(\R^n)\otimes S^p(T^*_{x_i}M)\big)$.
Since  $\lb D\a,\a\rb\geq 0$ and $B$ is a bundle map,  (\ref{E:(C)}) implies that
as $t\to\infty$, the eigenvalues (multiplicities counted) of the operator
$t^{-1}D_t$
concentrate near the eigenvalues of the model operator $K$ (see Proposition~1.2 of \cite{Sh}).
As a corollary, we obtain (P2) in the Introduction  as follows.

\begin{prop}[Corollary~1.3 of \cite{Sh}]\label{P:Main}
There exists $t_0>0$ such that for any $t>t_0$, we have
$$
\dim\ker D_t=\dim \ker t^{-1}D_t \leq \dim \ker K=\sum_i\dim\ker \bK(x_i).
$$
\end{prop}

\section{Hermite polynomials and operator formalism}

\label{HP-MO}

Throughout the remaining sections, let $x$ denote a critical point of the Morse function $f$
with Morse index $m$ as in Section~\ref{Adapted metric}. For notational simplicity,
we set
$$
\bK=\bK(x)\ \ \ \ \text{and}\ \ \ \ S^p=\Sym^p(T^*_xM).
$$
As the model operator $\bK$ acts in $L^2(\R^n)\otimes S^p$, Hermite functions on $\R^n$ are useful for calculating $\ker\bK$ because they constitute an orthogonal basis for $L^2(\R^n)$ (cf. \cite{F}).
However, we will use Hermite polynomials on $\R^n$ and relevant weighted operators for computational convenience.
In this section, we write weighted operators in terms of operators convenient for handling Hermite polynomials and use this operator formalism to prove Theorem~\ref{T:Main}\,(a).

Recall that for a  multi-index $I=(i_1,\cdots,i_n)$,
$$
I!=i_1!\cdots i_n!\ \ \ \ \text{and}\ \ \ \ |I|={\textstyle \sum} i_k.
$$
Also recall that $\{e^J\,|\,e^k=dy_k, |J|=p\}$ is an orthogonal basis for $S^p$.
Let $1_k$ denote the multi-index with $i$-th component $\d_{ik}$. Note that $I=\sum i_k 1_k$.

\subsection{Hermite polynomials on $\R^n$}

The classical {\em Hermite polynomials on $\R$} are given by
$$
H_k(t)=(-1)^ke^{t^2}\sdfrac{d^k}{dt^k}e^{-t^2}.
$$
Here the degree of $H_k(t)$ is $k$.
These Hermite polynomials on $\R$ satisfy
\begin{equation}\label{E:P-HP}
H_k^\p(t) = 2k H_{k-1}\ \ \ \ \text{and}\ \ \ \ tH_k(t)=kH_{k-1}(t)+\sdfrac12 H_{k+1}.
\end{equation}

\begin{defn}
\mbox{}
\begin{itemize}
\item[(a)]
The Hermite polynomials on $\R^n$ are products of Hermite polynomials on $\R$, that is,
$$
H_I(y)=H_{i_1}(y_1)\cdots H_{i_n}(y_n),
$$
where $I\in\Z^n$. If some $i_k<0$, we set $H_I=0$.
\item[(b)]
Denote the space generated by degree $q$ Hermite polynomials by
$$
\H^q=\sp\{H_I:|I|=q\}.
$$
\end{itemize}
\end{defn}

Below are well-known properties of the Hermite polynomials $H_I$ that are frequently used for our subsequent discussions.
\begin{itemize}

\item
For $d\mu=e^{-|y|^2}dy$,  $\{H_I\}$ is an orthogonal basis for
$L^2(\R^n,d\mu)$; hence the vectors $H_Ie^J$ form an orthogonal basis
for $L^2(\R^n,d\mu)\otimes S^p$ such that
\begin{equation}\label{E:Int}
\big( H_Ie^J,H_{I^\p}e^{J^\p}\big)=
\int_{\R^n} H_IH_{I^\p}d\mu \,\lb e^J,e^{J^\p}\rb  = \d_{I,I^\p}\d_{J,J^\p}\big(\sqrt{\pi}\big)^n 2^{|I|}I! J!.
\end{equation}

\item
The dimension of the space $\H^q$ is the number of multi-indices $I$ with $|I|=q$, so $\dim \H^q=n_q$, where $n_q$ is given in (\ref{E:dim}).

\item
By (\ref{E:P-HP}), the Hermite polynomials $H_I$
satisfy
\begin{equation}\label{E:key1}
\pa_k H_I=2i_k H_{I-1_k}\ \ \ \ \text{and}\ \ \ \ y_k H_I=i_kH_{I-1_k}+\sdfrac12 H_{I+1_k}.
\end{equation}
\end{itemize}

\subsection{Weighted operators}

Noting that $\phi\in L^2(\R^n,d\mu)$ if and only if $e^{-\frac12|y|^2}\phi\in L^2(\R^n)$, we define:

\begin{defn}
For the exponential weight function $w=e^{-\frac12|y|^2}$, let
$$
\bK_w=w^{-1}\bK w,\ \text{that is},\ K_w(\b)=w^{-1}\bK(w\b).
$$
\end{defn}

The weighted operator $\bK_w$ acts in  $L^2(\R^n,d\mu)\otimes S^p$ such that
$$
\dim\ker\bK_w=\dim\ker\bK\ \ \ \ \text{and}\ \ \ \
\bK_w=w^{-1}\bD w + \bB + \bV,
$$
where by a simple calculation
\begin{equation}\label{E:Dw}
w^{-1}\bD w
= -\sum_{k,\ell}(-\d_{k\ell}+y_k y_\ell -
y_k\pa_\ell - y_\ell\pa_k +\pa_k\pa_\ell  )(a_kc_\ell+c_k a_\ell).
\end{equation}

\subsection{Canonical commutation relations}

Without further confusion, we will often regard $y_k$ as an operator
$$
y_k:g(y)\mapsto y_kg(y).
$$
Noting (\ref{E:key1}), we obtain handy operators as follows.

\begin{defn}\label{D:CA}
For $k=1,\cdots, n$, set
$$
C_k=2y_k  -\pa_k\ \ \ \ \text{and}\ \ \ \
A_k=\sdfrac12\pa_k.
$$
\end{defn}

$C_k$ and $A_k$ are analogous to $c_k$ and $a_k$ in the following sense. Recall that
\begin{equation}\label{E:handy1}
\begin{aligned}
c_k:S^p&\to S^{p+1}:e^J\mapsto e^{J+1_k},
\\
a_k:S^p&\to S^{p-1}:e^J\mapsto j_ke^{J-1_k},
\end{aligned}
\end{equation}
Analogously, when restricted to $\H^q$, by (\ref{E:key1}) $C_k$ and $A_k$ are given as
\begin{equation}\label{E:handy2}
\begin{aligned}
&C_k:\H^q\to\H^{q+1}:H_I\mapsto H_{I+1_k},
\\
&A_k:\H^q\to \H^{q-1}:H_I\mapsto i_kH_{I-1_k}.
\end{aligned}
\end{equation}
Furthermore, one has:

\begin{lemma}\label{L:CCR-H}
\mbox{}
\begin{itemize}
\item[(a)]
$A_k^*=\frac12 C_k$ with respect to the inner product of $L^2(\R^n,d\mu)$.
\item[(b)]
The operators $C_k$ and $A_k$ satisfy the canonical commutation relations:
\begin{equation}\label{E:CCR-H}
C_kC_\ell=C_\ell C_k,\ \ \ \ A_kA_\ell=A_\ell A_k,\ \ \ \
C_k A_\ell -A_\ell C_k = -\d_{k\ell}I.
\end{equation}
\end{itemize}
\end{lemma}


The proof is immediate from definitions.

\subsection{Operator formalism}


\begin{defn}
\mbox{}
\begin{itemize}
\item[(a)]
Recalling $s_k$ in (\ref{E:simple-cal}) is defined by the Morse index $m$, we set
$$
(Cc)^-=\sum_{k\leq m}C_kc_k,\ \ \ \ (Cc)^+=\sum_{m<\ell}C_\ell c_\ell,
\ \ \ \ (Ca)^-=\sum_{k\leq m}C_ka_k,
$$
and likewise for $(Ca)^+, (Ac)^-, (Ac)^+, (Aa)^-, (Aa)^+$.
\item[(b)]
Let $P=(Ca)^--2(Aa)^+$ and $Q=(Cc)^--2(Ac)^+$.
\end{itemize}
\end{defn}

The proposition below is our key fact for analyzing $\ker \bK_w$. In particular, it shows $\bK_w$ is self-adjoint and positive semi-definite. Observe that as $A_k^*=\frac12 C_k$ and $a_k^*=c_k$, we have
$$
P^*=2(Ac)^--(Cc)^+\ \ \ \ \text{and}\ \ \ \ Q^*=2(Aa)^--(Ca)^+.
$$

\begin{prop}\label{P:WMO}
$\bK_w=P^*P+Q^*Q$ and
\begin{equation}\label{E:CR-Q}
Q^*Q-QQ^*=2\sum_k C_kA_k + 2\sum_{k\leq m}c_ka_k -2\sum_{m<\ell}c_\ell a_\ell + 2m.
\end{equation}
\end{prop}

\pf Since $\pa_k=2A_k$ and $y_k=A_k+\frac12 C_k$, by (\ref{E:CCR-H}) and (\ref{E:Dw}) one has
\begin{align}\label{A:Step1}
\bK_w &= \sum_{k,\ell\leq m}(C_kA_\ell+C_\ell A_k)(a_k c_\ell + c_k a_\ell)
+\sum_{m<k,\ell}(C_kA_\ell+C_\ell A_k)(a_k c_\ell + c_k a_\ell)
\notag\\
&+2\sum_{k\leq m}(a_kc_k + c_ka_k) -2\sum_{k\leq m<\ell}(4A_kA_\ell + C_kC_\ell)(a_kc_\ell+c_ka_\ell).
\end{align}
Since $C_\ell A_k(a_kc_\ell +c_k a_\ell)=C_\ell A_k(c_\ell a_k + a_\ell c_k)$ by (\ref{E:CCR}),
the equation (\ref{A:Step1}) becomes
\begin{align}\label{A:Step2}
\bK_w &= 2(Ca)^-(Ac)^-+2(Cc)^-(Aa)^-+2(Ca)^+(Ac)^++2(Cc)^+(Aa)^+
\notag\\[3pt]
&-4(Aa)^-(Ac)^+-(Cc)^-(Ca)^+-4(Ac)^-(Aa)^+-(Ca)^-(Cc)^+
\notag\\
&+2\sum_{k\leq m}(a_kc_k+c_ka_k).
\end{align}
On the other hand, by (\ref{E:CCR}) and (\ref{E:CCR-H}), one obtain
\begin{equation}\label{E:Step3}
\begin{aligned}
(Ca)^-(Ac)^-&=(Ac)^-(Ca)^-+\sum_{k\leq m}A_kC_k-\sum_{k\leq m}c_ka_k -m,
\\
(Cc)^-(Aa)^- &=(Aa)^-(Cc)^--\sum_{k\leq m}A_kC_k-\sum_{k\leq m}a_kc_k +m,
\\
(Ca)^-(Cc)^+&=(Cc)^+(Ca)^-,\ \ \ \ (Cc)^-(Ca)^+=(Ca)^+(Cc)^-.
\end{aligned}
\end{equation}
Now, plugging (\ref{E:Step3}) to the equation (\ref{A:Step2}) and then factoring gives the first assertion.
For the second assertion, apply (\ref{E:CCR}) and (\ref{E:CCR-H}) to have
$$
2(Ca)^+(Ac)^+=2(Ac)^+(Ca)^++2\sum_{m<\ell}A_\ell C_\ell -2\sum_{m<\ell}c_\ell a_\ell.
$$
This and the second equation of (\ref{E:Step3}) then yields the second assertion. \mqed

\begin{rem}\label{R:unweighted}
The weighted operator $\bK_w=w^{-1}\bK w$ is positive semi-definite with respect to the inner product on $L^2(\R^n,d\mu)\otimes S^p$, so the model operator
$\bK$ is positive semi-definite with respect to the inner product on $L^2(\R^n)\otimes S^p$.
\end{rem}

\medskip
As a corollary to Proposition~\ref{P:WMO}, we have:

\begin{cor}\label{C:m=n}
If $m=n$, then $\ker\bK_w=0$.
\end{cor}

\pf Let $m=n$. Then $(Aa)^+=(Ac)^+=0$, hence we have
$$
\bK_w=(Ca)^*Ca + (Cc)^*Cc,
$$
where  $Ca=\sum_kC_ka_k$ and $Cc=\sum_k C_k c_k$.
As $Cc$ is injective, we conclude $\bK_w$ is positive definite, hence $\ker\bK_w=0$.
\mqed

\subsection{Proof of Theorem~\ref{T:Main}\,(a)}

By Corollary~\ref{C:m=n} and $\dim\ker\bK=\dim\ker\bK_w$, the proposition below completes the proof of Theorem~\ref{T:Main}\,(a).

\begin{prop}
If $m\geq p$, then $\ker\bK_w=0$.
\end{prop}

\pf
Suppose $\b\in\ker \bK_w$. Then by Proposition~\ref{P:WMO}, $Q^*Q(\b)=0$. Noting  that $QQ^*$ is positive semi-definite, by (\ref{E:CR-Q}) we have
\begin{align}\label{A:ineq1}
0 &\geq \Big(\pig(\sum_k C_kA_k   + \sum_{k\leq m} c_ka_k -\sum_{m<\ell}c_\ell a_\ell +m\pig)\b,\b\Big)
\notag \\
&=2\sum_k\big(A_k\b,A_k\b) + \sum_{k\leq m}(a_k\b,a_k\b) - \sum_{m<\ell}\big(a_\ell\b,a_\ell\b) + m||\b||^2.
\end{align}
Since $\b$ can be written as $\b=\sum_{q,J}\b^q_J e^J$, where $\b^q_J e^J\in\H^q\otimes S^p$,
by (\ref{E:Int}), (\ref{E:handy1}) and (\ref{E:handy2}), the inequality (\ref{A:ineq1}) implies that
\begin{equation*}\label{E:ineq2}
0\geq \sum_{q,J} \Big(q+\sum_{k\leq m}j_k -\sum_{m<\ell} j_\ell+ m\Big)||\b^q_Je^J||^2.
\end{equation*}
As $\sum j_\ell\leq |J|=p\leq m$, it follows  that
$$
0\geq \sum_{q,J} \Big(q+\sum_{k\leq m}j_k \Big)||\b^q_Je^J||^2.
$$
This inequality implies that $\b^q_J=0$ for all $q\ne 0$, hence $\b=\sum_J\b^0_Je^J\in\H^0\otimes S^p$ such that
$$
0=Q(\b)=(Cc)^-\b-2(Ac)^+\b=(Cc)^-\b.
$$
Since $(Cc)^-$ is injective on $\H^0\otimes S^p$ (as $m\ne 0$), we conclude that $\b=0$. \mqed

\section{Dimension count}

This section proves Theorem~\ref{T:Main}\,(b) and (c) in the Introduction.

\subsection{Weighted operator with Morse index zero}

Let $\bK_w$ be the weighted  operator as in Section~\ref{HP-MO}. From now on, we assume that
the Morse index $m=0$. In this case, as $(Ca)^-=(Cc)^-=0$, by Proposition~\ref{P:WMO},
$$
\bK_w=(2Aa)^*2Aa + (2Ac)^*2Ac, 
$$
where $Aa=\sum_k A_ka_k$ and $Ac=\sum_k A_kc_k$. Consider
$$
\xymatrix{
\H^{q+1}\otimes S^{p-1} \ar@<0.5ex>[r]^(0.55){Ac}
&\H^q\otimes S^p \ar@<0.5ex>[r]^(0.45){Ac} \ar@<0.5ex>[l]^(0.45){Ca} \ar[d]^(0.45){Aa}
& \H^{q-1}\otimes S^{p+1}  \ar@<0.5ex>[l]^(0.55){Ca}
\\
& \H^{q-1}\otimes S^{p-1}
}
$$
and observe that:
\begin{itemize}
\item[(a)]
$\bK_w(\H^q\otimes S^p)\subset \H^q\otimes S^p$.
\item[(b)]
If $\b=\sum_q\b_q\in\ker\K_w$, where $\b_q\in\H^q\otimes S^p$, then by (a),
$$
\b_q\in\ker\bK_w,\ \ \ \forall q\geq 0.
$$
\item[(c)]
$\ker\bK_w=\ker(Aa)\cap \ker(Ac)$. In particular, $\H^0\otimes S^p\subset \ker \bK_w$.

\item[(d)]
By (\ref{E:handy1}) and (\ref{E:handy2}), it follows from (\ref{E:CR-Q}) that
\begin{equation}\label{E:EV}
(CaAc-AcCa)|_{\H^q\otimes S^p} = (q-p)Id.
\end{equation}
\item[(e)]
Since both $CaAc$ and $AcCa$ are positive semi-definite, (\ref{E:EV}) shows
\begin{equation}\label{E:qp}
q>p  \Rightarrow \ker Ac|_{\H^q\otimes S^p}=0
\ \ \ \text{and}\ \ \
q<p  \Rightarrow \ker Ca|_{\H^q\otimes S^p}=0.
\end{equation}

\end{itemize}

With this understood, we define:

\begin{defn}
Let $\bK_w$ be as above with $m=0$ and define
$$
\K^{q,p}=  
\ker (Aa|_{\H^q\otimes S^p})\cap \ker (Ac|_{\H^q\otimes S^p}).
$$
\end{defn}

One may then use the above facts (a) -- (e) to conclude:

\begin{lemma}\label{L:DC}
Let $\bK_w$ be as above with $m=0$. Then
$$
\ker\bK_w =\bigoplus_{q\leq p} \K^{q,p}\ \ \ \ \text{and}\ \ \ \ \K^{0,p}=\H^0\otimes S^p.
$$
\end{lemma}

\subsection{Proof of Theorem~\ref{T:Main}\,(b)}

Let $p=1$ and $m=0$. In this case, by Lemma~\ref{L:DC},
$$
K^{n,1}_0=\dim\ker\bK=\dim\ker\bK_w=\dim \H^0\otimes S^1+\dim\K^{1,1}.
$$
Since $\dim \H^0\otimes S^1=n$, Theorem~\ref{T:Main}\,(b) follows from the following lemma.

\begin{lemma}
${\ds \dim\K^{1,1}=\binom{n}{2}}$.
\end{lemma}

\pf
Write $\b\in  \H^1\otimes S^1$ as
$\b=\sum_{k,\ell}\b^k_\ell H_k e^\ell$. Then
$$
Ac(\b)
=\sum_{k,\ell}\b^k_\ell H_0 e^ke^\ell=0\ \ \Longleftrightarrow\ \ \b^k_\ell+\b^\ell_k=0,\ \forall k,\ell.
$$
This implies that $\ker Ac|_{\H^1\otimes S^1}$ is spanned by the vectors $H_k e^\ell-H_\ell e^k$. The lemma now follows since $Aa(H_k e^\ell-H_\ell e^k)=0$ for all $k\ne\ell$. \mqed.

\subsection{Proof of Theorem~\ref{T:Main}\,(c)}

For $p=2$ and $m=0$, by Lemma~\ref{L:DC},
$$
K^{n,2}_0=\dim\ker\bK=\dim\ker\bK_w=\dim \H^0\otimes S^2+\dim\K^{1,2}+\dim\K^{2,2}.
$$
Since ${\ds \dim \H^0\otimes S^2=\binom{n+1}{2} }$, Theorem~\ref{T:Main}\,(c) follows from Lemma~\ref{L:12} and
Lemma~\ref{L:22} below.


\begin{lemma}\label{L:12}
$\K^{1,2}=0$ for $n=1$, and if $n\geq 2$, we have
\begin{equation*}\label{E:dim-12}
\dim\K^{1,2}= n{n+1\choose 2}-{n+2\choose 3}-n.
\end{equation*}
\end{lemma}

\pf
Since the map
$Ca:\H^0\otimes S^3\to \H^1\otimes S^2$
is injective by (\ref{E:qp}), its adjoint map is surjective. Hence,
\begin{equation}\label{E:dim-count}
\dim\ker Ac|_{\H^1\otimes S^2}
= n{n+1\choose 2}-{n+2\choose 3}.
\end{equation}
If $n=1$, $\ker Ac|_{\H^1\otimes S^2}$ is trivial by (\ref{E:dim-count}), hence so is $\K^{1,2}$.
For $n\geq 2$, by the dimension formula (\ref{E:dim-count}), it suffices to show that
the restriction map $Aa:\ker Ac|_{\H^1\otimes S^2}\to \H^0\otimes S^1$ is surjective.
Let $u_{k\ell}=H_ke^ke^\ell-H_\ell e^ke^k$. Then
$$
Ac(u_{k\ell})=0\ \ \ \ \text{and}\ \ \ \ Aa(u_{k\ell})=e^\ell.
$$
This shows the restriction map $Aa:\ker Ac|_{\H^1\otimes S^2}\to \H^0\otimes S^1$ is surjective. \mqed


\begin{lemma}\label{L:22}
$\K^{2,2}=0$ for $n=1,2$, and if $n\geq 3$, we have
\begin{equation}\label{E:dim-22}
\dim\K^{2,2}={n+1 \choose 2}{n+1 \choose 2} - n{n+2 \choose 3} - {n+1 \choose 2}.
\end{equation}
\end{lemma}

\pf
For notational simplicity, let $Ac=Ac|_{\H^2\otimes S^2}$.
Since the map
$$
Ca:\H^1\otimes S^3\to \H^2\otimes S^2
$$
is injective by (\ref{E:qp}),  its adjoint $Ac$ is surjective; hence we have
\begin{equation*}\label{E:dim-qp22}
\dim \ker Ac
={n+1 \choose 2}{n+1 \choose 2} - n{n+2 \choose 3}.
\end{equation*}
This shows $\K^{2,2}=0$ when $n=1$. Let $n\geq 2$.
One may choose, by inspection,  linearly independent vectors in $\ker Ac$ as follows.
For distinct $i,j,k,\ell$, let
\begin{align*}
u_{k\ell}&=2H_{k\ell}e^{k\ell}-H_{\ell\ell}e^{kk}-H_{kk}e^{\ell\ell}\ (n\geq 2),
\\
v_{k\ell j}&=H_{kk}e^{\ell j}- H_{k\ell}e^{kj}-H_{kj}e^{k\ell}+H_{\ell j}e^{kk}\ (n\geq 3),
\\
w_{ijk\ell}&=H_{ij}e^{k\ell}-2H_{ik}e^{j\ell}+H_{i\ell}e^{jk}+H_{jk}e^{i\ell}-2H_{j\ell}e^{ik}+H_{k\ell}e^{ij}
\ (n\geq 4).
\end{align*}
Here $H_{k\ell}=H_{(k,\ell)}$ and $e^{k\ell}=e^{(k,\ell)}$. These vectors form a basis for $\ker Ac$ because they span a subspace of dimension
$$
{n\choose 2} + 3{n\choose 3} + 2{n\choose 4} =  {n+1 \choose 2}{n+1 \choose 2} - n{n+2 \choose 3} = \dim\ker Ac.
$$

Observing that
$$
Aa(u_{k\ell})=2H_\ell e^\ell+2H_k e^k,\ \
Aa(v_{k\ell j})=-H_\ell e^j-H_je^\ell,\ \ Aa(w_{ijk\ell})=0,
$$
we define a subspace of $\H^1\otimes S^1$ by
$$
S=\lb H_ke^k+H_\ell e^\ell\,|\,k\ne \ell\rb \oplus \lb H_ke^\ell+H_\ell e^k\,|\,k\ne \ell\rb\subset \H^1\otimes S^1.
$$
Obviously, $Aa(\ker Ac)\subset S$ and the restriction map
$Aa:\ker Ac\to S$ is injective for $n=2$, and surjective for $n\geq 3$.
Therefore, the lemma follows from counting dimension of $S$: for $n\geq 3$,
$$
\dim\lb H_ke^k+H_\ell e^\ell\,|\,k\ne \ell\rb =n\ \ \ \text{and}\ \ \
\dim\lb H_ke^\ell+H_\ell e^k\,|\,k\ne \ell\rb ={n \choose 2},
$$
where the first equality follows since $Aa(u_{k\ell}+u_{kj}-u_{\ell j})=2H_ke^k$. \mqed


\non
{\em Department of Mathematics, Dongguk University, Seoul 04620, Korea}

\non
{\em kwanghochoi\@@dongguk.edu}

\medskip

\non
{\em Department of  Mathematics,  University of Central Florida, Orlando, FL 32816}

\non
{\em Korea Institute for Advanced Study, Seoul 02455, Korea}

\non
{\em junho.lee\@@ucf.edu}

\end{document}